\documentclass[reqno,12pt]{amsart}
\usepackage{txfonts,mathrsfs}
\usepackage{cite}
\usepackage{amssymb}
\usepackage{amsfonts}
\usepackage{amsbsy}
\usepackage{graphicx}

\newtheorem{thm}{Theorem}

\newtheorem{prop}[thm]{Proposition}
\newtheorem{cor}[thm]{Corollary}

\newcommand{\Z}{{\mathbb Z}}

\def\tr{\mathop{\rm Tr}\nolimits}

\begin{document}

\title[gerneralized Thue-Morse Hamiltonians]
{The Hausdorff dimension of spectrum of a class of gerneralized Thue-Morse Hamiltonians}

\author{Qinghui Liu}
\address{Sch. Comp. Sci., Beijing Inst. Tech., Beijing 100081, China}
\email[Correponding author]{qhliu@bit.edu.cn}

\author{Zhiyi Tang}
\email{zytang@bit.edu.cn}

\begin{abstract}
We study a class of Schr\"odinger operators $H_{m,\lambda}$ with
generalized Thue-Morse potential that generated by the substitution
$\tau(a)=a^mb^m$, $\tau(b)=b^ma^m$ on two symbol alphabet $\Sigma=\{a,b\}$
for integer $m\ge 2$ and coupling $\lambda>0$.
We show that
$$\dim_H \sigma(H_{m,\lambda})\ge \frac{\log \Lambda_m}{\log 64m+4},$$
where $\sigma(H_{m,\lambda})$ is the spectrum of $H_{m,\lambda}$, $\Lambda_2=2$, and for $m>2$,
$\Lambda_m=m$, if $m\equiv0\mod 4$;
$\Lambda_m=m-3$, if $m\equiv1\mod 4$;
$\Lambda_m=m-2$, if $m\equiv2\mod 4$;
$\Lambda_m=m-1$, if $m\equiv3\mod 4$.
This implies that $\dim_H \sigma(H_{m,\lambda})$ tends to $1$ as $m$ tends to infinity.



\end{abstract}

\maketitle

\section{Introduction}

For many classes of Schr\"odinger operators with quasi-periodic potential
taking finitely values, the corresponding spectrum has zero Lesbegue measure,
see \cite{Su89,BIST89,BG93,L02,LTWW02,DL06,LQ11,LQ12}.
Among them, potentials generated by Sturmian sequences, substitution sequences are closely studied.

Let $H_{\alpha,\lambda}$ be the Schr\"odinger operator with Sturmian potential,
where $0<\alpha<1$ is irrational frequency and $\lambda>0$ is the coupling.
If $\alpha=(\sqrt{5}-1)/2$, i.e., the Fibonacci potential,
then the Hausdorff dimension of $\sigma(H_{\alpha,\lambda})$
tends to $0$ as the coupling $\lambda$ tends to infinity(see \cite{LW04,DEGT08,LPW07}),
and the Hausdorff dimension of $\sigma(H_{\alpha,\lambda})$
tends to $1$ as the coupling $\lambda$ tends to $0$(see \cite{DG11}).
Letting $\alpha=[0;a_1,a_2,\cdots]$ be its continued fractional expansion,
if $$\liminf_{n\rightarrow\infty}\sqrt[n]{a_1a_2\cdots a_{n}}=\infty,$$
then, for all large $\lambda$, $\dim_H\sigma(H_{\alpha,\lambda})=1$(see \cite{LW04}).

Let $H_{tm,\lambda}$ and $H_{pd,\lambda}$ be the Schr\"odinger operator
with corresponding Thue-Morse potential and period-doubling potential and
with the coupling $\lambda>0$. It is shown in \cite{LQ15, LQY21} that,
$$\dim_H\sigma(H_{tm,\lambda})>\frac{\ln 2}{140\ln 2.1},\
\dim_H\sigma(H_{pd,\lambda})>\frac{\ln (\sqrt{5}+1)/2}{\ln 4}.
$$
The authors get the results by studying the trace maps
generated by the Thue-Morse substitution and the period-doubling substitution respectively.
In \cite{LQ15}, the authors focus on looking for a set of germs.
In \cite{LQY21}, to get the lower bound of the Hausdorff dimension,
the authors focus on looking for a separating nested structure(SNS).
Both of the papers heavily rely on details of the dynamical properties of
the corresponding trace maps.
However, without studying details of the corresponding dynamical system,
we can find SNS for some generalized Thue-Morse potential.
And then, we can get a common positive lower bound
of Hausdorff dimension of the spectrum for any coupling.

In \cite{KA91}, Kolar and Ali studied a class of generalized Thue-Morse substitution
$$a\rightarrow a^mb^n,\ b\rightarrow b^n a^m,$$
for $m,n\ge1$, and the corresponding potential.
We focus on a special class of generalized Thue-Morse substitution, i.e., taking $m\ge1$,
$$\tau(a)=a^mb^m,\ \tau(b)=b^ma^m.$$
Note that the classical Thue-Morse substitution correspond to the case $m=1$.
It is seen that
$\tau^{2n}(a)$ is both a prefix  and suffix of $\tau^{2(n+1)}(a)$.
Define a two-sided sequence $\xi$ as
$$\xi:=\lim_{n\rightarrow\infty}\tau^{2n}(a)|\tau^{2n}(a)
=\cdots \xi(-2)\xi(-1)|\xi(0)\xi(1)\cdots.$$
Define the potential $V_\xi=(V_\xi(n))_{n\in\Z}$ by,
$V_\xi(n)=1$ if $\xi(n)=a$ and $V_\xi(n)=-1$ if $\xi(n)=b$ for $n\in\mathbb{Z}$.
Take $\lambda>0$.
Let $H_{m,\lambda}$ be the discrete Schr\"odinger operator acting
 on $\ell^2(\Z)$ with potential $\lambda V_\xi$, i.e.,
for any $n\in\mathbb{Z}$,
$$(H_{m,\lambda}\psi)_n=\psi_{n+1}+\psi_{n-1}+\lambda V_\xi(n)\psi_n,
\ \ \ \forall \psi\in\ell^2(\Z).$$
Denote by $\sigma(H_{m,\lambda})$ the spectrum of $H_{m,\lambda}$.

Define $\Lambda_2=2$ and, for $m>2$, define
$\Lambda_m=m$, if $m\equiv0\mod 4$;
$\Lambda_m=m-3$, if $m\equiv1\mod 4$;
$\Lambda_m=m-2$, if $m\equiv2\mod 4$;
$\Lambda_m=m-1$, if $m\equiv3\mod 4$.
Define $\gamma_2=8(5+\sqrt{29})$ and, for $m>2$, define
$\gamma_m=1+32m+\sqrt{1+192m+1024m^2}$.

\begin{thm}\label{hdim-m}
For any integer $m\ge2$ and $\lambda>0$,
$$\dim_H \sigma(H_{m,\lambda})\ge \frac{\log \Lambda_m}{\log \gamma_m}
> \frac{\log \Lambda_m}{\log 64m+4}.$$
\end{thm}

\noindent
{\bf Remark}. The method in this paper cannot apply directly to the case of $m=1$.
By detail analysis on dynamical property of the Thue-Morse trace maps as \cite{LQY21},
we can also improve the lower bound given in \cite{LQ15}.

\begin{cor}\label{hdim-infty}
For any $\lambda>0$,
$$\lim_{m\rightarrow\infty}\dim_H \sigma(H_{m,\lambda})=1.$$
\end{cor}

The rest of the paper are organized as following.
In Section 2, we show a trace map of $\tau$ given by Kolar and Ali in \cite{KA91},
and prove some properties of Chebyshev polynomial needed in our proof.
In Section 3, we prove Theorem \ref{hdim-m} in case of $m>2$.
In Section 4, we prove Theorem \ref{hdim-m} in case of $m=2$.

\section{Trace polynomials and Chebyshev polynomials}
Fix $\lambda>0$. The transfer matrix on site $n\in\mathbb{Z}$ is often defined as
$$T_n:=\left[
\begin{array}{cc}
	t-\lambda V_\xi(n)&-1\\
	1&0
\end{array}\right].$$
We have $\tr T_n\cdots T_2T_1$ is a polynomial on $t$ with degree $n$,
where $\tr A$ be trace of a matrix $A$.
These polynomials are called trace polynomials.
For $n>0$, define
\begin{equation}\label{sptrace}
x_n:=\tr T_{(2m)^n}\cdots T_2T_1,\ \sigma_n:=\{t\in\mathbb{R}\ |\ |x_n(t)|\le2\}.
\end{equation}
By \cite{LTWW02} (see also \cite{BG93,L02}),
\begin{equation}\label{sp}
\sigma(H_{m,\lambda})=\overline{\bigcap_{n\ge1}\bigcup_{k\ge n} \sigma_k},
\end{equation}
where $\overline{A}$ is the closure of a subset $A\subset\mathbb{R}$.

In general, it is difficult to fix $\sigma(H_{m,\lambda})$.
For Thue-Morse potential, i.e., the case $m=1$,
\begin{equation}\label{1in2}
\sigma_n\cup\sigma_{n+1}\supset \sigma_{n+2},\ \forall n\ge1.\end{equation}
Then $$\sigma(H_{1,\lambda})=\bigcap_{n\ge1}\sigma_n\cup\sigma_{n+1}.$$
We do not know whether \eqref{1in2} hold for other generalized Thue-Morse substitutions.
But for the substitution $a\rightarrow abb$, $b\rightarrow bba$, \eqref{1in2} do not hold.
We will construct a subset of $\sigma(H_{m,\lambda})$ by \eqref{sp},
and then estimate the lower bound of Hausdorff dimension of the subset.

\begin{figure}[htbp]
\centering
\includegraphics[width=0.5\textwidth]{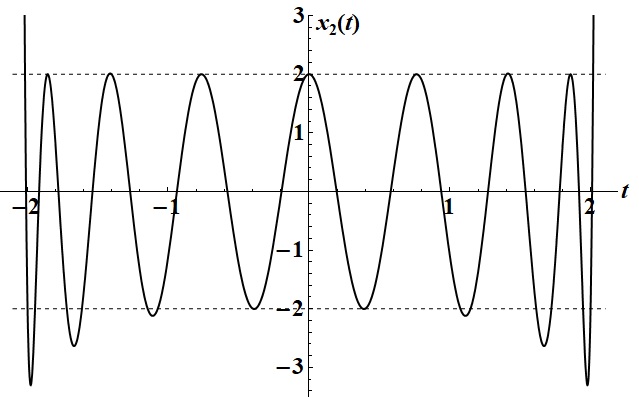}
\caption{The graph of $x_2(t)$ in case of $m=2$ and $\lambda=0.1$.}\label{fig1}
\end{figure}

Let $p(t)$ be a trace polynomial on $t$ with degree $k$. By Floquet theory or \cite{T89},
$p(t)=2$ has $k$ real roots, and $p(t)=-2$ also has $k$ real roots.
By \cite{T89}, if $|p(t)|<2$ then $p'(t)\ne0$. This implies
$\{t\ |\ |p(t)|\le2\}$ is composed of $k$ non-overlapping intervals.
Moreover, on each of these intervals $I$,  $p(I)=[-2,2]$ and $p(t)$ is monotone.
As a corollary,
\begin{equation}\label{Toda}
p''(t)<0\quad \mbox{ if } p'(t)=0 \mbox{ and } p(t)\ge2.
\end{equation}
Figue \ref{fig1} shows the graph of $x_2(t)$ in case of $m=2$ and $\lambda=0.1$.

To study $(x_n)_{n\ge1}$, we redefine them in the following dynamical system.
Define matrix functions on real variables $t$,
$$
M_0:= \left[
\begin{array}{cc}
	t-\lambda&-1\\
	1&0
\end{array}
\right]\ \ \ \text{ and }\ \ \ N_0:= \left[
\begin{array}{cc}
	t+\lambda&-1\\
	1&0
\end{array}
\right].
$$
For any $n\ge0$, let
$$M_{n+1}=N_n^m M_n^m,\ \ N_{n+1}=M_n^m N_n^m.$$
 Define, for any $n\ge0$,
$$x_n:=\tr M_n,\ \ y_n:=\tr M_nN_n.$$
It is clear that, for $n\ge1$, $\tr M_n=\tr N_n$. Fix $\lambda>0$.
We denote by $x_n', x_n'',y_n'$
the first or second derivation of $x_n$, $y_n$ on $t$.

By equation (19) in \cite{KA91}, for any $n\ge1$,
\begin{equation}\label{ka19}
\left\{
\begin{array}{l}
x_{n+1}=(d_m(x_n))^2(y_n-2)+2\\
y_{n+1}=(d_{2m}(x_n))^2(y_n-2)+2,
\end{array}\right.
\end{equation}
where $d_k$ are the Chebyshev polynomials of the second kind,
i.e., $d_0(t)\equiv0$, $d_1(t)\equiv1$,
and for $k\ge1$,
$$d_{k+1}(t)=t\ d_k(t)-d_{k-1}(t).$$

For $k>0$, $\theta\in[0,\pi]$,
\begin{equation}\label{cheb<=2}
d_k(2\cos\theta)=\frac{\sin k\theta}{\sin\theta}.
\end{equation}
This implies that for $1\le j\le k-1$,
$$d_k\left(2\cos\frac{j\pi}{k}\right)=0.$$

For any $t\in[2\cos3\pi/4,2\cos\pi/4]=[-\sqrt{2},\sqrt{2}]$,
there is $\theta\in[\pi/4,3\pi/4]$ such that $2\cos\theta=t$.
Then \eqref{cheb<=2} implies
\begin{equation}\label{dkupb}
|d_k(t)|=\frac{|\sin k\theta|}{|\sin\theta|} \le \frac{1}{\sqrt{2}/2}=\sqrt{2}.\end{equation}
Taking derivation on \eqref{cheb<=2}, we can get
$$d_k'(2\cos\theta)=\frac{d}{d\theta}\left(\frac{\sin k\theta}{\sin\theta}\right)\frac{-1}{2\sin\theta}=
\frac{\sin(k-1)\theta-(k-1)\cos k\theta\sin\theta}{2\sin^3\theta}.$$
Then, for $t\in[-\sqrt{2},\sqrt{2}]$,
\begin{equation}\label{ddkupb}
|d_k'(t)|\le\frac{k}{2(\sqrt{2}/2)^3}= \sqrt{2}k.\end{equation}

Since
$$\frac{d_{2k}(2\cos\theta)}{d_k(2\cos\theta)}=\frac{\sin 2k\theta}{\sin k\theta}=2\cos k\theta,$$
we have for $t\in[-2,2]$,
\begin{equation}\label{chebyprop}
 |d_{2k}(t)|\le 2|d_k(t)|.
\end{equation}

\section{Proof of Theorem \ref{hdim-m} in case of $m>2$}
For $n\ge1$, define $\mathscr{F}_n$ be the set of all intervals $[u,v]$
such that $x_n$ is monotone on $[u,v]$ and
$$x_n([u,v])=[-2,2].$$
Since $x_n(t)$ is a trace polynomial of degree $(2m)^n$,
$\sharp \mathscr{F}_n=(2m)^n$.

In the following, we construct a seperating nested structure(SNS) $(\mathscr{G}_n)_{n\ge1}$ in $(\mathscr{F}_n)_{n\ge1}$.
That is, for any $n\ge1$, $\mathscr{G}_n\subset\mathscr{F}_n$;
for any intervals $I,J\in\mathscr{G}_n$, $I\cap J$ is empty or contain a single point;
for any $J\in\mathscr{G}_{n+1}$, there is a unique interval $I\in\mathscr{G}_n$ so that $J\subset I$.
The limit set of $(\mathscr{G}_n)_{n\ge1}$ is defined as
$$\bigcap_{n\ge1}\bigcup_{I\in\mathscr{G}_n}I.$$

For $n=1$, we take an interval $I\in\mathscr{F}_1$, and set $\mathscr{G}_1=\{I\}$.

Suppose, for $n\ge1$, $\mathscr{G}_n$ is defined.
Take an interval $I=[u,v]\in\mathscr{G}_n$.
Then $x_n$ is monotone on $[u,v]$ and
$x_n([u,v])=[-2,2].$
For $1\le j\le m-1$, there is unique $t_j\in[u,v]$ such that $x_n(t_j)=2\cos j\pi/m$.
By \eqref{ka19},
$$x_{n+1}(t_j)=2,\ x_{n+1}'(t_j)=0.$$
By \eqref{Toda}, $x_{n+1}''(t_j)<0$. Hence, there are $t_{j-1}<u_j,v_j<t_{j+1}$ such that
$$x_{n+1}(u_j)=x_{n+1}(v_j)=-2,$$
and $x_{n+1}(t)$ is increasing on $(u_j,t_j)$ and decreasing on $(t_j,v_j)$.
It is clear that the intervals $[u_j,t_j]$, $[t_j,v_j]\in\mathscr{F}_{n+1}$.
For $2\le j\le m-2$, $$[u_j,t_j], [t_j,v_j]\subset [u,v].$$
And also, $[t_1,v_1]$, $[u_{m-1},t_{m-1}]\subset [u,v]$.

Note that, for $t\in[u_j,v_j]$, $-2\le x_{n+1}(t)\le 2$.
By \eqref{ka19},
$$-4\le d_m^2(x_n(t))(y_n(t)-2)\le0.$$
Then, by \eqref{chebyprop}, $|d_{2m}^2(x_n(t))(y_n(t)-2)|\le 16$, i.e.,
$$|y_{n+1}(t)-2|\le16.$$
Let $m_0=\lceil m/4\rceil$, $m_1=\lfloor 3m/4\rfloor$.
Define $$P_I=P_{[u,v]}=\{[t_{m_0},v_{m_0}], [u_{m_1},t_{m_1}]\}\cup
\bigcup_{j=m_0+1}^{m_1-1}\{[u_j,t_j], [t_j,v_j]\}\subset\mathscr{F}_{n+1}.$$
It is direct that, for $m>2$,
\begin{equation}\label{count}
\sharp P_{[u,v]}=\Lambda_m.\end{equation}

Note that, for $t\in \bigcup_{J\in P_I}J$,
letting $x_n(t)=2\cos\theta$, we have $\pi/4<\theta<3\pi/4$.
Hence, by \eqref{dkupb} and \eqref{ddkupb},
$$|d_m(x_n(t))|\le\sqrt{2},\ |d_m'(x_n(t))|\le \sqrt{2}m,\  |d_{2m}'(x_n(t))|\le 2\sqrt{2}m.$$

Define
$$\mathscr{G}_{n+1}=\bigcup_{I\in\mathscr{G}_n}P_I.$$
Continue this process, we get a SNS $\{\mathscr{G}_n\}_{n\ge1}$.
In the SNS, each interval in $\mathscr{G}_n$ contain
$\Lambda_m$ intervals in $\mathscr{G}_{n+1}$.

As a summary, for any $J=[u,v]\in\mathscr{G}_{n}$,
for $1<k<n$, $t\in[u,v]$,
$$\begin{array}{l}
|x_k(t)|\le2,\ |y_k(t)-2|\le16,\\
|d_m(x_k(t))|\le\sqrt{2},\ |d_m'(x_k(t))|\le \sqrt{2}m,\\
|d_{2m}(x_k(t))|\le\sqrt{2},\ |d_{2m}'(x_k(t))|\le 2\sqrt{2}m.
\end{array}$$
By \eqref{ka19}, for $t\in[u,v]$ and $1\le k<n$,
$$\begin{array}{rcl}
x_{k+1}'(t)&=&2d_m(x_k(t))d_m'(x_k(t))(y_k(t)-2)x_k'(t)+d_m^2(x_k(t))y_k'(t),\\
y_{k+1}'(t)&=&2d_{2m}(x_k(t))d_{2m}'(x_k(t))(y_k(t)-2)x_k'(t)+d_{2m}^2(x_k(t))y_k'(t).
\end{array}$$
Then, we have,
$$\begin{array}{rcl}
|x_{k+1}'(t)|&\le& 64m|x_k'(t)|+2|y_k'(t)|,\\
|y_{k+1}'(t)|&\le& 128m|x_k'(t)|+2|y_k'(t)|.
\end{array}$$
For $m>2$, let $\gamma_m=1+32m+\sqrt{1+192m+1024m^2}(<64m+4)$,
which is the largest eigenvalue of the coefficient matrix
$$\left[\begin{array}{cc}64m&2\\ 128m&2\end{array}\right].$$
Then, there exist $M>0$ independent of $n$ and $J$ such that, for $t\in J$,
$$|x_{n}'(t)|\le M \gamma_m^n.$$
Then
\begin{equation}\label{intlen}
4=\left|\int_u^v x_n'(t)dt\right|=\int_u^v |x_n'(t)|dt\le M\gamma_m^n(v-u).
\end{equation}
And hence
\begin{equation}\label{len}
|J|=v-u\ge 4M^{-1}\gamma_m^{-n}.\end{equation}

Let $E$ be the limit set of the SNS.
By our construction, for each $n\ge1$ and each $I\in\mathscr{G}_n$, $I\subset\sigma_n$.
Then by \eqref{sp}, $E\subset \sigma(H_{m,\lambda})$.

The following is the Proposition 5.7 in \cite{LQY21}.

\begin{prop}\label{SNS-lower-dim}
Let $(\mathscr{I}_n)_{n\ge1}$ be a SNS. Assume it satisfies:

i) There exist $w\in (0,1)$ and $C>0$ such that
$$
|I|\ge Cw^n,\ \ \forall n\ge0, \forall I\in \mathscr{I}_n.
$$

ii) There exists $C'\ge1$ such that for any $n,k\ge0$, and  $I,I'\in \mathscr{I}_n$,
$$
\frac{\#\{J\in \mathscr{I}_{n+k}: J\subset I\}}{\#\{J\in \mathscr{I}_{n+k}: J\subset I'\}}\le C'.
$$

Let $A$ be the limit set of $(\mathscr{I}_n)_{n\ge1}$. Then
\begin{equation*}
  \dim_H A\ge \liminf_{n\to\infty}\frac{\log \# \mathscr{I}_n}{-n\log w}.
\end{equation*}
\end{prop}

It is clear that, for the SNS $(\mathscr{G}_n)_{n\ge1}$,
by \eqref{len}, $w=\gamma_m^{-1}$;
by \eqref{count}, $C'=1$ and $\sharp\mathscr{G}_n=\Lambda_m^{n-1}$.
Then
$$\dim_H E\ge \frac{\log\Lambda_m}{\log \gamma_m}>\frac{\log\Lambda_m}{\log 64m+4}.$$

\section{Proof of Theorem \ref{hdim-m} in case of $m=2$}
In this case, since $d_2(t)=t$, $d_4(t)=t^3-2t$, by \eqref{ka19}, for $n\ge1$,
\begin{equation}\label{m=2}
\left\{
\begin{array}{l}
x_{n+1}=(d_2(x_n))^2(y_n-2)+2=x_n^2(y_n-2)+2\\
y_{n+1}=(d_{4}(x_n))^2(y_n-2)+2=(x_n^3-2x_n)^2(y_n-2)+2.
\end{array}\right.
\end{equation}
It is direct that
\begin{equation*}
y_{n+1}-2=(x_{n+1}-2)(x_n^2-2)^2.
\end{equation*}
Then, for any $n\ge1$,
\begin{equation}\label{reccurence}
x_{n+2}=x_{n+1}^2(x_{n+1}-2)(x_{n}^2-2)^2+2.
\end{equation}

It is direct that, for $n\ge2$, if $x_n(t)=0$ or $x_n(t)=2$ for some $t$, then
$$x_{n+1}(t)=2,\ x_{n+1}'(t)=0,\ x_{n+1}''(t)<0.$$

For $n\ge1$, define again $\mathscr{F}_n$ be the set of all intervals $[u,v]$
such that $x_n$ is monotone on $[u,v]$ and
$$x_n([u,v])=[-2,2].$$
It is clear that $\sharp \mathscr{F}_n=4^n$.
We construct another seperating nested structure(SNS)
$(\mathscr{G}_n)_{n\ge1}$ in $(\mathscr{F}_n)_{n\ge1}$.

For $n=1$, we take an interval $I\in\mathscr{F}_1$, and set $\mathscr{G}_1=\{I\}$.

Suppose, for $n\ge1$, $\mathscr{G}_n$ is defined.
Take an interval $J=[u,v]\in\mathscr{G}_n$.
Then $x_n$ is monotone on $[u,v]$ and
$x_n([u,v])=[-2,2].$
Take $u\le t_0<t_1\le v$ so that $x_n(t_i)=0$ or $2$ for $i=0,1$.
Since $x_{n+1}(t_0)=x_{n+1}(t_1)=2$, there exist $t_0<\hat{t}_0\le \hat{t}_1<t_1$ such that
$x_{n+1}(\hat{t}_0)=x_{n+1}(\hat{t}_1)=-2$ and
$$\{[t_0,\hat{t}_0],[\hat{t}_1,t_1]\}\subset \mathscr{F}_{n+1}.$$
Let $P_J=\{[t_0,\hat{t}_0],[\hat{t}_1,t_1]\}$. Define
$$\mathscr{G}_{n+1}=\bigcup_{I\in\mathscr{G}_n}P_I.$$
It is clear that $\sharp \mathscr{G}_{n+1}=2^n$.
We get a SNS $(\mathscr{G}_n)_{n\ge1}$. Let $E$ be the limit set of the SNS.
Then $E\subset\sigma(H_{2,\lambda})$.

By \eqref{reccurence},
$$x_{n+2}'= x_{n+1}(3x_{n+1}-4)(x_n^2-2)^2 x_{n+1}'+4x_nx_{n+1}^2(x_{n+1}-2)(x_n^2-2)x_n'.$$

For $n\ge1$, take any interval $J=[u,v]\in \mathscr{G}_n$.
Since, by the construction above, for $1\le k\le n$ and $t\in [u,v]$, $|x_k(t)|\le2$,
we have $|3x_{k}(t)-4|\le10$, $|x_k(t)-2|\le4$, $|x_k^2(t)-2|\le2$.
This implies for $1< k< n$ and $t\in [u,v]$,
$$|x_{k+1}'(t)|\le 80|x_{k}'(t)|+256|x_{k-1}'(t)|.$$
Let $\gamma_2=8(5+\sqrt{29})$, which is the largest root of $x^2-80x-256=0$.
Then there exists $M$ independent of $n$ and $J$ so that, for $t\in J$,
$$|x_n'(t)|\le M \gamma_2^n.$$
By a similar argument as \eqref{intlen},
$$|J|=v-u\ge 4M^{-1}\gamma_2^{-n}.$$

The SNS $(\mathscr{G}_n)_{n\ge1}$ satisfy condition of Proposition \ref{SNS-lower-dim}
with parameters $w=\gamma_2^{-1}$, $C'=1$ and $\sharp\mathscr{G}_n=\Lambda_2^{n-1}$.
Then
$$\dim_H \sigma(H_{2,\lambda})\ge\dim_H E\ge \frac{\log\Lambda_2}{\log \gamma_2}
>\frac{\log 2}{\log 88}.$$


\medskip

\noindent
{\bf Acknowledgement}. The authors thank Prof. Wen Zhiying and Qu Yanhui
for helpful discussions. The authors are supported by
National Natural Science Foundation of China No.11871098.


\begin{thebibliography}{99}

\bibitem{BIST89}
J. Bellissard, B. Iochum, E. Scoppola, D. Testard, Spectral properties of one dimensional
quasi-crystals, Commun. Math. Phys. 125:3(1989), 527-543.

\bibitem{BG93}
 A. Bovier, J. Ghez, Spectral properties of one-dimensional Schr\"odinger operators with potentials generated by substitutions, Commun. Math. Phys. 158:1(1993), 45--66.

\bibitem{DEGT08}
D. Damanik, M. Embree, A. Gorodetski, S. Tcheremchantsev, The fractal dimension of the spectrum
of the Fibonacci Hamiltonian, Comm. Math. Phys., 280:2(2008), 499--516.

\bibitem{DG11}
D. Damanik, A. Gorodetski,
Spectral and quantum dynamical properties of the weakly coupled Fibonacci Hamiltonian,
Commun. Math. Phys. 305:1(2011), 221-277.

\bibitem{DL06}
D. Damanik, D. Lenz, A condition of Boshernitzan and uniform convergence
in the multiplicative ergodic theorem, Duke Math. J., 133:1(2006), 95--123.

\bibitem{L02}
D. Lenz, Singular spectrum of Lebesgue measure zero for one-dimensional quasicrystals,
Commun. Math. Phys. 227:1(2002), 119--130.

\bibitem{LQ11}
Q.-H. Liu, Y.-H. Qu, Uniform convergence of Schr\"odinger cocycles over simple
Toeplitz subshift, Annales Henri Poincar\'e 12:1(2011), 153--172.

\bibitem{LQ12}
Q.-H. Liu, Y.-H. Qu, Uniform convergence of Schr\"odinger cocycles over bounded
Toeplitz subshift, Annales Henri Poincar\'e 13:6(2012), 1483--1500.

\bibitem{LQ15}  Q.-H. Liu,  Y.-H. Qu,
On the Hausdorff dimension of the spectrum of Thue-Morse Hamiltonian,
Comm. Math. Phys., 338:2(2015), 867-891.

\bibitem{LQY21}  Q.-H. Liu,  Y.-H. Qu, X. Yao,
The spectrum of period-doubling Hamiltonian, preprint, arXiv:2108.13257.

\bibitem{LPW07}
Q.H. Liu, J. Peyri\`ere, Z.Y. Wen,
Dimension of the spectrum of one-dimensional discrete Schr\"odinger
operators with Sturmian potentials, C. R. Math., 345:12(2007), 667--672.

\bibitem{LTWW02}
Q. H. Liu, B. Tan, Z. X. Wen, J. Wu,
Measure zero spectrum of a class of Schr\"odinger operators,
J. Stat. Phys. 106:3-4(2002), 681-691.

\bibitem{LW04}
Q. H. Liu, Z. Y. Wen,
Hausdorff dimension of spectrum of one-dimensional Schr\"{o}dinger operator with Sturmian potentials,
Potential Analysis, 20:1(2004), 33-59.

\bibitem{KA91}
M. Kolar, M. K. Ali,
Generalized Thue-Morse chains and their physical properties,
Physical Review B, 43:1(1991), 1034-1047.

\bibitem{Su89}
A. S\"uto, Singular continuous spectrum on a Cantor set of zero Lebesgue
measure for the Fibonacci hamiltonian, J. Stat. Phys. 56:3-4(1989), 525-531.

\bibitem{T89}
M. Toda, Theory of Nonlinear Lattices, 2nd enlarged edn, Solid-State Sciences 20, Springer-
Verlag, 1989, Chap. 4.

\end{thebibliography}
\end{document}